\documentclass[12pt,leqno]{article}
\usepackage{amssymb,graphics,latexsym}
\usepackage{graphicx}
\usepackage{amsmath}
\large\normalsize

\def\R{\hbox{\bf\rlap{I}{\hbox to 2 pt{}}R}}
\def\tb#1#2{\mathop{#1\vphantom{\sum}}\limits_{\displaystyle #2}}

\newcommand{\re}{{\rm Re\, }}

\newcommand{\dia}{{\rm diag\, }}

\newcommand{\rank}{{\rm rank\, }}
\newcommand{\nul}{{\rm nullity\, }}

\begin{document}
\thispagestyle{empty}
\begin{center}
\section*{Structures and Numerical Ranges of Power Partial Isometries}

\vspace*{3mm}

\begin {tabular}{lcl}
\hspace*{1cm}{\bf Hwa-Long Gau}$^{*1}$\hspace*{1cm}&and & \hspace*{1cm}{\bf Pei Yuan Wu}$^2$
\vspace*{3mm}\\
Department of Mathematics & & Department of Applied Mathematics\\
National Central University&& National Chiao Tung University\\
Chung-Li 32001, Taiwan&& Hsinchu 30010, Taiwan\\
Republic of China&&Republic of China
\end{tabular}

\end{center}

\vspace{3mm}

\centerline{\bf Abstract}

We derive a matrix model, under unitary similarity, of an $n$-by-$n$ matrix $A$ such that $A, A^2, \ldots, A^k$ ($k\ge 1$) are all partial isometries, which generalizes the known fact that if $A$ is a partial isometry, then it is unitarily similar to a matrix of the form ${\scriptsize\left[\begin{array}{cc} 0 & B\\ 0 & C\end{array}\right]}$ with $B^*B+C^*C=I$. Using this model, we show that if $A$ has ascent $k$ and $A, A^2, \ldots, A^{k-1}$ are partial isometries, then the numerical range $W(A)$ of $A$ is a circular disc centered at the origin if and only if $A$ is unitarily similar to a direct sum of Jordan blocks whose largest size is $k$. As an application, this yields that, for any $S_n$-matrix $A$, $W(A)$ (resp., $W(A\otimes A)$) is a circular disc centered at the origin if and only if $A$ is unitarily similar to the Jordan block $J_n$. Finally, examples are given to show that the conditions that $W(A)$ and $W(A\otimes A)$ are circular discs at 0 are independent of each other for a general matrix $A$.

\vspace{3mm}

\noindent
\emph{AMS classification}: 15A99, 15A60\\
\emph{Keywords}: Power partial isometry, numerical range, $S_n$-matrix.

\vspace{3mm}

${}^*$Corresponding author.

E-mail addresses: hlgau@math.ncu.edu.tw (H.-L. Gau), pywu@math.nctu.edu.tw (P. Y. Wu)

${}^1$Research supported by the National Science Council of the Republic of China under NSC-102-2115-M-008-007.

${}^2$Research supported by the National Science Council of the Republic of China under NSC-102-2115-M-009-007 and by the MOE-ATU project.

\newpage
\noindent
{\bf\large 1. Introduction}

\vspace{4mm}

An $n$-by-$n$ complex matrix $A$ is a \emph{partial isometry} if $\|Ax\|=\|x\|$ for any vector $x$ in the orthogonal complement $(\ker A)^{\perp}$ in $\mathbb{C}^n$ of the kernel of $A$, where $\|\cdot\|$ denotes the standard norm in $\mathbb{C}^n$. The study of such matrices or, more generally, such operators on a Hilbert space dates back to 1962 \cite{6}. Their general properties have since been summarized in \cite[Chapter 15]{5}.

\vspace{4mm}

In this paper, we study matrices $A$ such that, for some $k\ge 1$, the powers $A, A^2, \ldots, A^k$ are all partial isometries. In Section 2 below, we derive matrix models, under unitary similarity, of such a matrix (Theorems 2.2 and 2.4). They are generalizations of the known fact that $A$ is a partial isometry if and only if it is unitarily similar to a matrix of the form ${\scriptsize\left[\begin{array}{cc} 0 & B\\ 0 & C\end{array}\right]}$ with $B^*B+C^*C=I$ (Lemma 2.1).

\vspace{4mm}

Recall that the \emph{ascent} of a matrix, denoted by $a(A)$, is the minimal integer $k\ge 0$ for which $\ker A^k=\ker A^{k+1}$. It is easily seen that $a(A)$ is equal to the size of the largest Jordan block associated with the eigenvalue 0 in the Jordan form of $A$. We denote the $n$-by-$n$ \emph{Jordan block}
$$\left[
    \begin{array}{cccc}
      0 & 1 &   &   \\
        & 0 & \ddots &   \\
        &   &  \ddots & 1 \\
        &   &   & 0
    \end{array}
  \right]$$
by $J_n$. The \emph{numerical range} $W(A)$ of $A$ is the subset $\{\langle Ax, x\rangle : x\in \mathbb{C}^n, \|x\|=1\}$ of the complex plane $\mathbb{C}$, where $\langle\cdot, \cdot\rangle$ is the standard inner product in $\mathbb{C}^n$. It is known that $W(A)$ is a nonempty compact convex subset, and $W(J_n)=\{z\in\mathbb{C} : |z|\le\cos(\pi/(n+1))\}$ (cf. \cite[Proposition 1]{4}). For other properties of the numerical range, the readers may consult \cite[Chapter 22]{5} or \cite[Chapter 1]{9}.

\vspace{4mm}

Using the matrix model for power partial isometries, we show that if $a(A)=k\ge 2$ and $A, A^2, \ldots, A^{k-1}$ are all partial isometries, then the following are equivalent: (a) $W(A)$ is a circular disc centered at the origin, (b) $A$ is unitarily similar to a direct sum $J_{k_1}\oplus J_{k_2}\oplus\cdots\oplus J_{k_{\ell}}$ with $k=k_1\ge k_2\ge\cdots\ge k_{\ell}\ge 1$, and (c) $A$ has no unitary part and $A^j$ is a partial isometry for all $j\ge 1$ (Theorem 2.6). An example is given, which shows that the number ``$k-1$'' in the above assumption is sharp (Example 2.7).

\vspace{4mm}

In Section 3, we consider the class of $S_n$-matrices. Recall that an $n$-by-$n$ matrix $A$ is of {\em class} $S_n$ if $A$ is a contraction ($\|A\|\equiv\max\{\|Ax\|: x\in\mathbb{C}^n, \|x\|=1\}\le 1$), its eigenvalues are all in $\mathbb{D}$ ($\equiv\{z\in\mathbb{C} : |z|<1\}$), and it satisfies $\rank(I_n-A^*A)=1$. Such matrices are the finite-dimensional versions of the \emph{compression of the shift} $S(\phi)$, first studied by Sarason \cite{10}. They also feature prominently in the Sz.-Nagy--Foia\c{s} contraction theory \cite{11}. It turns out that a hitherto unnoticed property of such matrices is that if $A$ is of class $S_n$ and $k$ is its ascent, then $A, A^2, \ldots, A^k$ are all partial isometries. Thus the structure theorems in Section 2 are applicable to $A$ or even to $A\otimes A$, the tensor product of $A$ with itself. As a consequence, we obtain that, for an $S_n$-matrix $A$, the numerical range $W(A)$ (resp., $W(A\otimes A)$) is a circular disc centered at the origin if and only if $A$ is unitarily similar to the Jordan block $J_n$ (Theorem 3.3). The assertion concerning $W(A)$ is known before (cf. \cite[Lemma 5]{12}). Finally, we give examples to show that if $A$ is a general matrix, then the conditions for the circularity (at the origin) of $W(A)$ and $W(A\otimes A)$ are independent of each other (Examples 3.5 and 3.6).

\vspace{4mm}

We use $I_n$ and $0_n$ to denote the $n$-by-$n$ identity and zero matrices, respectively. An identity or zero matrix with unspecified size is simply denoted by $I$ or $0$. For an $n$-by-$n$ matrix $A$, $\nul A$ is used for $\dim\ker A$, and $\rank A$ for its rank. The \emph{real part} of $A$ is $\re A=(A+A^*)/2$. The \emph{geometric} and \emph{algebraic multiplicities} of an eigenvalue $\lambda$ of $A$ are $\nul(A-\lambda I_n)$ and the multiplicity of the zero $\lambda$ in the characteristic polynomial $\det(zI_n-A)$ of $A$, respectively. An $n$-by-$n$ diagonal matrix with diagonal entries $a_1, \ldots, a_n$ is denoted by $\dia(a_1, \ldots, a_n)$.

\vspace{8mm}

\noindent
{\bf \large 2. Power Partial Isometries}

\vspace{4mm}

We start with the following characterizations of partial isometries.

\vspace{4mm}

{\bf Lemma 2.1.} \emph{The following conditions are equivalent for an $n$-by-$n$ matrix $A$}:

 (a) \emph{$A$ is a partial isometry},

 (b) \emph{$A^*A$ is an} (\emph{orthogonal}) \emph{projection}, \emph{and}

 (c) \emph{$A$ is unitarily similar to a matrix of the form ${\scriptsize\left[\begin{array}{cc} 0 & B\\ 0 & C\end{array}\right]}$ with $B^*B+C^*C=I$}.

\noindent
\emph{In this case}, ${\scriptsize\left[\begin{array}{cc} 0 & B\\ 0 & C\end{array}\right]}$ \emph{acts on} $C^n=\ker A\oplus(\ker A)^{\perp}$.

\vspace{4mm}

Its easy proof is left to the readers.

\vspace{4mm}

The next theorem gives the matrix model, under unitary similarity, of a matrix $A$ with $A, A^2, \ldots, A^k$ ($1\le k\le a(A)$) partial isometries.

\vspace{4mm}

{\bf Theorem 2.2.} \emph{Let $A$ be an $n$-by-$n$ matrix}, $\ell\ge 1$, \emph{and} $k=\min\{\ell, a(A)\}$. \emph{Then the following conditions are equivalent}:

(a) $A, A^2, \ldots, A^k$ \emph{are partial isometries},

(b) \emph{$A$ is unitarily similar to a matrix of the form}
$$A'\equiv\left[\begin{array}{ccccc} 0 & A_1 & & & \\ & 0 & \ddots & & \\ & & \ddots & A_{k-1} & \\ & & & 0 & B\\ & & & & C\end{array}\right] \ on \  \mathbb{C}^n=\mathbb{C}^{n_1}\oplus\cdots\oplus\mathbb{C}^{n_k}\oplus\mathbb{C}^{m},$$
\emph{where the} $A_j$'\emph{s satisfy $A_j^*A_j=I_{n_{j+1}}$ for $1\le j\le k-1$}, \emph{and $B$ and $C$ satisfy $B^*B+C^*C=I_m$}. \emph{In this case}, $n_j=\nul A$ \emph{if} $j=1$, $\nul A^j-\nul A^{j-1}$ \emph{if} $2\le j\le k$, \emph{and} $m=\rank A^k$,

(c) \emph{$A$ is unitarily similar to a matrix of the form}
$$A''\equiv\left[\begin{array}{ccccc} 0 & I & & & \\ & 0 & \ddots & & \\ & & \ddots & I & \\ & & & 0 & B\\ & & & & C\end{array}\right]\oplus(J_{k-1}\oplus\cdots\oplus J_{k-1})\oplus\cdots\oplus(J_1\oplus\cdots\oplus J_1)$$
$$on \  \mathbb{C}^n=\underbrace{\mathbb{C}^{n_k}\oplus\cdots\oplus\mathbb{C}^{n_k}}_{k}\oplus\mathbb{C}^{m}
\oplus\underbrace{\mathbb{C}^{k-1}\oplus\cdots\oplus\mathbb{C}^{k-1}}_{n_{k-1}-n_k}\oplus\cdots\oplus
\underbrace{\mathbb{C}\oplus\cdots\oplus\mathbb{C}}_{n_1-n_2},$$
\emph{where the} $n_j$'\emph{s}, $1\le j\le k$, \emph{and $m$ are as in} (b), \emph{and $B$ and $C$ satisfy $B^*B+C^*C=I_m$}.

\vspace{4mm}

For the proof of Theorem 2.2, we need the following lemma.

\vspace{4mm}

{\bf Lemma 2.3.} \emph{Let $A=[A_{ij}]_{i,j=1}^n$ be a block matrix with $\|A\|\le 1$}, \emph{and let $\alpha$ be a nonempty subset of} $\{1, 2, \ldots, n\}$. \emph{If for some} $j_0$, $1\le j_0\le n$, \emph{we have} $\sum_{i\in\alpha}A_{i j_0}^*A_{i j_0}=I$, \emph{then $A_{i j_0}=0$ for all $i$ not in $\alpha$}.

\vspace{4mm}

{\em Proof}. Since $\|A\|\le 1$, we have $A^*A\le I$. Thus the same is true for the $(j_0, j_0)$-block of $A^*A$, that is, $\sum_{i=1}^n A_{i j_0}^*A_{i j_0}\le I$. Together with our assumption that $\sum_{i\in\alpha}A_{i j_0}^*A_{i j_0}=I$, this yields $\sum_{i\not\in\alpha}A_{i j_0}^*A_{i j_0}\le 0$. It follows immediately that $A_{i j_0}=0$ for all $i$ not in $\alpha$.  \hspace{2mm} $\blacksquare$

\vspace{4mm}

{\em Proof of Theorem $2.2$}. To prove (a) $\Rightarrow$ (b), let $H_1=\ker A$, $H_j=\ker A^j\ominus\ker A^{j-1}$ for $2\le j\le\ell$, and $H_{\ell+1}=\mathbb{C}^n\ominus\ker A^{\ell}$. Note that if $\ell>a(A)$, then at most $H_1, \ldots, H_{k+1}$ are present. Hence $A$ is unitarily similar to the block matrix $A'\equiv[A_{ij}]_{i,j=1}^{k+1}$ on $\mathbb{C}^n=H_1\oplus\cdots\oplus H_{k+1}$. It is easily seen that $A_{ij}=0$ for any $(i,j)\neq(k+1,k+1)$ with $1\le j\le i\le k+1$. For the brevity of notation, let $A_j=A_{j, j+1}$, $1\le j\le k-1$, $B=A_{k, k+1}$, and $C=A_{k+1, k+1}$. We now check, by induction on $j$, that $A_j^*A_j=I_{n_{j+1}}$ for all $j$, and $A_{ij}=0$ for $1\le i\le j-2\le k-2$.

\vspace{4mm}

For $j=1$, since $A$ is a partial isometry, $A^*A$ is an (orthogonal) projection by Lemma 2.1. We obviously have $A^*A=0$ on $H_1=\ker A$ and $A^*A=I$ on $H_1^{\perp}=H_2\oplus\cdots\oplus H_{k+1}$. Thus $A'^*A'=0\oplus I\oplus\cdots\oplus I$ on $\mathbb{C}^n=H_1\oplus H_2\oplus\cdots\oplus H_{k+1}$. Since $A'^*A'$ is of the form
$$\left[
    \begin{array}{ccccc}
      0 & 0 & 0 & \cdots & 0 \\
      0 & A_1^*A_1 & * & \cdots & * \\
      0 & * & * & \cdots & * \\
      \vdots & \vdots & \vdots &   & \vdots \\
      0 & * & * & \cdots & *
    \end{array}
  \right],$$
we conclude that $A_1^*A_1=I$.

\vspace{4mm}

Next assume that, for some $p$ ($2\le p<k$), $A_j^*A_j=I$ for all $j$, $1\le j\le p-1$, and all the blocks in $A'$ which are above $A_1, \ldots, A_{p-1}$ are zero. We now check that $A_{p}^*A_{p}=I$ and all blocks above $A_{p}$ are zero. Since $A^{p}$ is a partial isometry, ${A^p}^*A^{p}$ is an (orthogonal) projection with kernel equal to $H_1\oplus\cdots\oplus H_{p}$. Thus ${A'^p}^*A'^{p}=\underbrace{0\oplus\cdots\oplus 0}_{p}\oplus \underbrace{I\oplus\cdots\oplus I}_{k-p+1}$. But from
$$A'=\left[
       \begin{array}{ccccccccc}
         0 & A_1 & 0 & \cdots & 0 & * & \cdots & * & * \\
           & \ddots & \ddots & \ddots & \vdots & \vdots &   & \vdots & \vdots \\
           &   & \ddots & \ddots & 0 & \vdots &   & \vdots & \vdots \\
           &   &   & \ddots & A_{p-1} & * &   & \vdots & \vdots \\
           &   &   &   & 0 & A_{p} & \ddots & \vdots & \vdots \\
           &   &   &   &   & 0 & \ddots & * & \vdots \\
           &   &   &   &   &   & \ddots & A_{k-1} & * \\
           &   &   &   &   &   &   & 0 & B \\
           &   &   &   &   &   &   &   & C \\
       \end{array}
     \right],$$
we have
$$A'^{p}=\begin{array}{ll} \ \ \ \overbrace{\ \hspace{15mm} \ }^{\displaystyle p} \ \ \ \overbrace{\ \hspace{71mm} \ }^{\displaystyle k-p+1} & \\ \left[
       \begin{array}{cccccccc}
         0 & \cdots & 0 & \prod_{j=1}^{p}A_j & * & \cdots & *  & * \\
         \cdot  &   &   & 0 & \prod_{j=2}^{p+1}A_j & \ddots & \vdots  &  \vdots \\
         \cdot  &   &   &   & \ddots  & \ddots & *  &  \vdots \\
         \cdot  &   &   &   &   & \ddots  & \prod_{j=k-p}^{k-1}A_j &  * \\
         \cdot  &   &   &   &   &   & 0 &  * \\
         \cdot  &   &   &   &   &   & \vdots &  \vdots \\
         \cdot  &   &   &   &   &   & 0 &  BC^{p-1} \\
         0 & \cdot  & \cdot  & \cdot  & \cdot & \cdot  & 0 &  C^{p} \\
       \end{array}
     \right] & \hspace{-11mm}\begin{array}{l}\left.\begin{array}{l} {\ } \\ {\ } \\ {\ }\\ {\ }\end{array}\right\}k-p\\ \left.\begin{array}{l}{\ } \\ {\ } \\ {\ }\\ {\ }\end{array}\right\}p+1\end{array}\end{array}.$$
Thus the $(p+1, p+1)$-block of ${A'^p}^*A'^{p}$ is $(\prod_{j=1}^{p}A_j)^*(\prod_{j=1}^{p}A_j)=A_{p}^*A_{p}$, which is equal to $I$ from above. Lemma 2.3 then implies that all the blocks in $A'$ which are above $A_{p}$ are zero. Thus, by induction, the first $k$ block columns of $A'$ are of the asserted form.

\vspace{4mm}

Finally, we check that $B^*B+C^*C=I_m$. If this is the case, then all the blocks in $A'$ above $B$ and $C$ are zero by Lemma 2.3 again and we will be done. As above, $A'^{k-1}$ is of the form
$$\left[\begin{array}{ccccc}
0 & \cdots & 0 & \prod_{j=1}^{k-1}A_j & D_1\\
0 & \cdots & 0 & 0 & D_2\\
\vdots &   & \vdots & \vdots & \vdots\\
0 & \cdots & 0 & 0 & D_k\\
0 & \cdots & 0 & 0 & C^{k-1}
\end{array}\right],$$
and the (orthogonal) projection ${A'^{k-1}}^*A'^{k-1}$ equals $\underbrace{0\oplus\cdots\oplus 0}_{k-1}\oplus I\oplus I$ on $\mathbb{C}^n=H_1\oplus\cdots\oplus H_{k-1}\oplus H_k\oplus H_{k+1}$. Hence the $(k+1, k+1)$-block of ${A'^{k-1}}^*A'^{k-1}$ is
\begin{equation}\label{e1}
(\sum_{j=1}^k D_j^*D_j)+{C^{k-1}}^*C^{k-1},
\end{equation}
which is equal to $I$. Similarly,
$$A'^k=A'^{k-1}A'=\left[\begin{array}{cccc}
0 & \cdots & 0 &   (\prod_{j=1}^{k-1}A_j)B+D_1C\\
0 & \cdots & 0 &   D_2C\\
\vdots &   & \vdots   & \vdots\\
0 & \cdots & 0 &   D_kC\\
0 & \cdots & 0 &   C^{k}
\end{array}\right]$$
and the $(k+1, k+1)$-block of ${A'^k}^*A'^k$,
\begin{equation}\label{e2}
B^*(\prod_{j=1}^{k-1}A_j)^*(\prod_{j=1}^{k-1}A_j)B+
B^*(\prod_{j=1}^{k-1}A_j)^*D_1C+C^*D_1^*(\prod_{j=1}^{k-1}A_j)B+(\sum_{j=1}^kC^*D_j^*D_jC)+{C^k}^*C^k,
\end{equation}
is also equal to $I_m$. We deduce from (\ref{e1}), (\ref{e2}) and $A_j^*A_j=I$ for $1\le j\le k-1$ that
\begin{equation}\label{e3}
B^*B+B^*(\prod_{j=1}^{k-1}A_j)^*D_1C+C^*D_1^*(\prod_{j=1}^{k-1}A_j)B+{C}^*C=I_m.
\end{equation}
To complete the proof, we need only show that $(\prod_{j=1}^{k-1}A_j)^*D_1=0$. Indeed, since $(\prod_{j=1}^{k-1}A_j)^*(\prod_{j=1}^{k-1}A_j)=I_{n_k}$, there is an $n_1$-by-$n_1$ unitary matrix $U$ such that $U^*(\prod_{j=1}^{k-1}A_j)={\scriptsize\left[\begin{array}{c} I_{n_k} \\ 0\end{array}\right]}$. Then $V\equiv U\oplus \underbrace{I\oplus\cdots\oplus I}_k$ is unitary and
$$V^*A'^{k-1}V=\left[\begin{array}{ccccc}
0 & \cdots & 0 & U^*(\prod_{j=1}^{k-1}A_j) & U^*D_1\\
0 & \cdots & 0 & 0 & D_2\\
\vdots &   & \vdots & \vdots & \vdots\\
0 & \cdots & 0 & 0 & D_k\\
0 & \cdots & 0 & 0 & C^{k-1}
\end{array}\right]=\left[\begin{array}{ccccc}
0 & \cdots & 0 & \left[\begin{array}{c} I_{n_k} \\ 0\end{array}\right] & \left[\begin{array}{c} 0 \\ D'_1\end{array}\right]\\
0 & \cdots & 0 & 0 & D_2\\
\vdots &   & \vdots & \vdots & \vdots\\
0 & \cdots & 0 & 0 & D_k\\
0 & \cdots & 0 & 0 & C^{k-1}
\end{array}\right].$$
Hence
$$(\prod_{j=1}^{k-1}A_j)^*D_1=\left[I_{n_k} \ 0\right]U^*U\left[\begin{array}{c} 0 \\ D'_1\end{array}\right]=\left[I_{n_k} \ 0\right]\left[\begin{array}{c} 0 \\ D'_1\end{array}\right]=0$$
as asserted. We conclude from (\ref{e3}) that $B^*B+C^*C=I_m$. Moreover, the sizes of the blocks in $A'$ are as asserted from our construction. This proves (a) $\Rightarrow$ (b).

\vspace{4mm}

Next we prove (b) $\Rightarrow$ (c). Let $A'$ be as in (b), and let $n_1, \ldots, n_k, m$ be the sizes of the diagonal blocks of $A'$. Since $A_j^*A_j=I_{n_{j+1}}$ for all $j$, $1\le j\le k-1$, we have $n_1\ge n_2\ge\cdots\ge n_k$. Also, from $A_{k-1}^*A_{k-1}=I_{n_k}$, we deduce that there is a unitary matrix $U_{k-1}$ of size $n_{k-1}$ such that $U_{k-1}^*A_{k-1}={\scriptsize\left[\begin{array}{c} I_{n_k} \\ 0\end{array}\right]}$. Similarly, since $(A_{k-2}U_{k-1})^*(A_{k-2}U_{k-1})=I_{n_{k-1}}$, there is a unitary $U_{k-2}$ of size $n_{k-2}$ such that $U_{k-2}^*(A_{k-2}U_{k-1})={\scriptsize\left[\begin{array}{c} I_{n_{k-1}} \\ 0\end{array}\right]}$. Proceeding inductively, we obtain a unitary $U_j$ of size $n_j$ satisfying $U^*_j(A_jU_{j+1})={\scriptsize\left[\begin{array}{c} I_{n_{j+1}} \\ 0\end{array}\right]}$ for each $j$, $1\le j\le k-3$. If $U=U_1\oplus\cdots\oplus U_{k-1}\oplus I_{n_k}\oplus I_m$, then
$$U^*A'U=\left[\begin{array}{cccccc} 0 & U_1^*A_1U_2 & & & & \\ & 0 & \ddots & & & \\ & & \ddots & U_{k-2}^*A_{k-2}U_{k-1} & & \\ & & & 0 & U_{k-1}^*A_{k-1} & \\ & & & & 0 & B\\ & & & & & C\end{array}\right]$$
$$=\left[\begin{array}{cccccc} 0 & \left[\begin{array}{c} I_{n_{2}} \\ 0\end{array}\right] & & & & \\ & 0 & \ddots & & & \\ & & \ddots & \left[\begin{array}{c} I_{n_{k-1}} \\ 0\end{array}\right] & & \\ & & & 0 & \left[\begin{array}{c} I_{n_{k}} \\ 0\end{array}\right] & \\ & & & & 0 & B\\ & & & & & C\end{array}\right].$$
Note that this last matrix is unitarily similar to the one asserted in (c).

\vspace{4mm}

To prove (c) $\Rightarrow$ (a), we may assume that
$$A''=\left[
        \begin{array}{ccccc}
          0 & I &   &   &   \\
            & 0 & \ddots &   &   \\
            &   & \ddots & I &   \\
            &   &   & 0 & B \\
            &   &   &   & C
        \end{array}
      \right]$$
with $B^*B+C^*C=I_m$. This is because powers of any Jordan block are all partial isometries and the direct sums of partial isometries are again partial isometries. Simple computations show that
$$A''^j=\begin{array}{ll} \ \ \ \overbrace{\ \hspace{15mm} \ }^{\displaystyle j} \  \overbrace{\ \hspace{36mm} \ }^{\displaystyle k-j+1} & \\ \left[
       \begin{array}{cccccccc}
         0 & \cdots & 0 & I & 0 & \cdots & 0  & 0 \\
           & \cdot  &   & 0 & \ddots & \ddots & \vdots  &  \vdots \\
           &   & \cdot  &   & \ddots  & \ddots & 0  &  \vdots \\
           &   &   & \cdot  &   & \ddots  & I &  0 \\
           &   &   &   & \cdot  &   & 0 &  B \\
           &   &   &   &   & \cdot  & \vdots &  \vdots \\
           &   &   &   &   &   & 0 &  BC^{j-1} \\
           &   &   &   &   &   &   &  C^{j} \\
       \end{array}
     \right] & \hspace{-11mm}\begin{array}{l}\left.\begin{array}{l} {\ } \\ {\ } \\ {\ }\\ {\ }\end{array}\right\}k-j\\ \left.\begin{array}{l}{\ } \\ {\ } \\ {\ }\\ {\ }\end{array}\right\}j+1\end{array}\end{array}$$
and ${A''^j}^*A''^j=\underbrace{0\oplus\cdots\oplus 0}_j\oplus \underbrace{I\oplus\cdots\oplus I}_{k-j}\oplus D$, where $D=(\sum_{s=0}^{j-1}{C^s}^*B^*BC^s)+{C^j}^*C^j$ for each $j$, $1\le j\le k$. From $B^*B+C^*C=I_m$, we deduce that
\begin{align*}
D &= B^*B+(\sum_{s=1}^{j-2}{C^s}^*B^*BC^s)+{C^{j-1}}^*(B^*B+C^*C)C^{j-1}\\
&= B^*B+(\sum_{s=1}^{j-2}{C^s}^*B^*BC^s)+{C^{j-1}}^*C^{j-1}\\
&= B^*B+(\sum_{s=1}^{j-3}{C^s}^*B^*BC^s)+{C^{j-2}}^*(B^*B+C^*C)C^{j-2}\\
&= \cdots\\
&= B^*B+C^*C\\
&= I_m.
\end{align*}
Hence ${A''^j}^*A''^j=0\oplus I$, which implies that $A''^j$ is a partial isometry by Lemma 2.1 for all $j$, $1\le j\le k$. This proves (c) $\Rightarrow$ (a). \hspace{2mm} $\blacksquare$

\vspace{4mm}

A consequence of Theorem 2.2 is the following.

\vspace{4mm}

{\bf Theorem 2.4.} \emph{Let $A$ be an $n$-by-$n$ matrix and $\ell>a(A)$}. \emph{Then the following conditions are equivalent}:

(a) $A, A^2, \ldots, A^{\ell}$ \emph{are partial isometries},

(b) \emph{$A$ is unitarily similar to a matrix of the form} $U\oplus J_{k_1}\oplus\cdots\oplus J_{k_m}$, \emph{where $U$ is unitary and} $a(A)=k_1\ge \cdots\ge k_m\ge 1$, \emph{and}

(c) \emph{$A^j$ is a partial isometry for all $j\ge 1$}.

\vspace{4mm}

The equivalence of (b) and (c) here is the finite-dimensional version of a result of Halmos and Wallen \cite[Theorem]{7}.

\vspace{4mm}

{\em Proof of Theorem $2.4$}. Since $\ell>k\equiv a(A)$, Theorem 2.2 (a) $\Rightarrow$ (b) says that $A$ is unitarily similar to a matrix of the form
$$A'\equiv\left[\begin{array}{ccccc} 0_{n_1} & A_1 & & & \\ & 0_{n_2} & \ddots & & \\ & & \ddots & A_{k-1} & \\ & & & 0_{n_k} & B\\ & & & & C\end{array}\right] \  \  \ \mbox{on} \  \ \mathbb{C}^n=\mathbb{C}^{n_1}\oplus\cdots\oplus\mathbb{C}^{n_k}\oplus\mathbb{C}^{m}$$
with the $A_j$'s, $B$ and $C$ satisfying the properties asserted therein. As $k$ is the ascent of $A$, $\nul A^k$ equals the algebraic multiplicity of eigenvalue 0 of $A$. Since $\nul A^k=\nul A'^k=\sum_{j=1}^kn_j$, it is seen from the structure of $A'$ that the eigenvalue 0 appears fully in the diagonal $0_{n_j}$'s. This shows that 0 cannot be an eigenvalue of $C$ or $C$ is invertible.

\vspace{4mm}

A simple computation yields that
$$A'^{k+1}=\left[\begin{array}{cccc}
0 & \cdots & 0 &   (\prod_{j=1}^{k-1}A_j)BC\\
0 & \cdots & 0 &   (\prod_{j=2}^{k-1}A_j)BC^2\\
\vdots &   & \vdots   & \vdots\\
0 & \cdots & 0 &   A_{k-1}BC^{k-1}\\
0 & \cdots & 0 &   BC^k\\
0 & \cdots & 0 &   C^{k+1}
\end{array}\right]$$
and
\begin{equation}\label{e4}
{A'^{k+1}}^*A'^{k+1}=0_{n_1}\oplus\cdots\oplus 0_{n_k}\oplus D,
\end{equation}
where, after simplification by using $A_j^*A_j=I_{n_{j+1}}$ for $1\le j\le k-1$, $D=(\sum_{j=1}^k{C^j}^*B^*BC^j)+{C^{k+1}}^*C^{k+1}$. As $A'^{k+1}$ is a partial isometry, ${A'^{k+1}}^*A'^{k+1}$ is a projection by Lemma 2.1. Moreover, we also have
$$\nul {A'^{k+1}}^*A'^{k+1}=\nul A'^{k+1}=\nul A'^k=\sum_{j=1}^kn_j,$$
where the second equality holds because of $k=a(A')$. Thus we obtain from (\ref{e4}) that $D=I_m$. Therefore,
\begin{align*}
I_m &= D = (\sum_{j=1}^k{C^j}^*B^*BC^j)+{C^{k+1}}^*C^{k+1}\\
&= (\sum_{j=1}^{k-1}{C^j}^*B^*BC^j)+{C^{k}}^*(B^*B+C^*C)C^{k}\\
&= (\sum_{j=1}^{k-1}{C^j}^*B^*BC^j)+{C^{k}}^*C^{k}\\
&= \cdots\\
&= C^*(B^*B+C^*C)C\\
&= C^*C.
\end{align*}
This shows that $C$ is unitary and hence $B=0$ (from $B^*B+C^*C=I_m$). Thus $A'$ is unitarily similar to the asserted form in (b). This completes the proof of (a) $\Rightarrow$ (b). The implications (b) $\Rightarrow$ (c) and (c) $\Rightarrow$ (a) are trivial. \hspace{2mm} $\blacksquare$

\vspace{4mm}

At this juncture, it seems appropriate to define the {\em power partial isometry index} $p(\cdot)$ for any matrix $A$:
$$p(A)\equiv\sup\{k\ge 0: I, A, A^2, \ldots, A^k \ \mbox{are all partial isometries}\}.$$
An easy corollary of Theorem 2.4 is the following estimate for $p(A)$.

\vspace{4mm}

{\bf Corollary 2.5.} \emph{If $A$ is an $n$-by-$n$ matrix}, \emph{then $0\le p(A)\le a(A)$ or $p(A)=\infty$}. \emph{In particular}, \emph{we have} (a) $0\le p(A)\le n-1$ \emph{or} $p(A)=\infty$, \emph{and} (b) $p(A)=n-1$ \emph{if and only if $A$ is unitarily similar to a matrix of the form}
\begin{equation}\label{e5}
\left[\begin{array}{ccccc}
0 & 1 &   &  &\\
  & 0 & \ddots &   &\\
  &   & \ddots & 1 &\\
  &   &   &  0 & a  \\
  &   & & &  b
\end{array}\right]\end{equation}
\emph{with $|a|^2+|b|^2=1$ and $a, b\neq 0$}.

\vspace{4mm}

{\em Proof}. The first assertion follows from Theorem 2.4. If $p(A)=n$, then $a(A)=n$, which implies that the Jordan form of $A$ is $J_n$. Thus $p(A)=\infty$, a contradiction. This proves (a) of the second assertion.

\vspace{4mm}

As for (b), if $p(A)=n-1$, then $a(A)=n$ will lead to a contradiction as above. Thus we must have $a(A)=n-1$. Theorem 2.2 implies that $A$ is unitarily similar to a matrix of the form (\ref{e5}) with $|a|^2+|b|^2=1$. Since either $a=0$ or $b=0$ will lead to the contradicting $p(A)=\infty$, we have thus proven one direction of (b). The converse follows easily from Theorem 2.2 and the arguments in the preceding paragraph. \hspace{2mm} $\blacksquare$

\vspace{4mm}

The next theorem gives conditions for which $p(A)\ge a(A)-1$ implies that $A$ is unitarily similar to a direct sum of Jordan blocks.

\vspace{4mm}

{\bf Theorem 2.6.} \emph{Let $A$ be an $n$-by-$n$ matrix with $p(A)\ge a(A)-1$}. \emph{Then the following conditions are equivalent}:

(a) \emph{$W(A)$ is a circular disc centered at the origin},

(b) \emph{$A$ is unitarily similar to a direct sum of Jordan blocks},

(c) \emph{$A$ has no unitary part and $A^j$ is a partial isometry for all $j\ge 1$}, \emph{and}

(d) \emph{$A$ has no unitary part and $A, A^2, \ldots, A^{\ell}$ are partial isometries for some $\ell>a(A)$}.

\noindent
\emph{In this case}, $W(A)=\{z\in \mathbb{C} : |z|\le\cos(\pi/(a(A)+1))\}$ \emph{and} $p(A)=\infty$.

\vspace{4mm}

Here a matrix is said to have {\em no unitary part} if it is not unitarily similar to one with a unitary summand.

\vspace{4mm}

Note that, in the preceding theorem, the condition $p(A)\ge a(A)-1$ cannot be replaced by the weaker $p(A)\ge a(A)-2$. This is seen by the next example.

\vspace{4mm}

{\bf Example 2.7.} If $A=J_3\oplus{\scriptsize\left[\begin{array}{cc} 0 & (1-|\lambda|^2)^{1/2}\\ 0 & \lambda\end{array}\right]}$, where $0<|\lambda|\le\sqrt{2}-1$, then $a(A)=3$ and $W(A)=\{z\in\mathbb{C} : |z|\le\sqrt{2}/2\}$. Since $A$ is a partial isometry while $A^2$ is not, we have $p(A)=1$. Note that $A$ has a nonzero eigenvalue. Hence it is not unitarily similar to any direct sum of Jordan blocks.

\vspace{4mm}

The proof of Theorem 2.6 depends on the following series of lemmas, the first of which is a generalization of \cite[Theorem 1]{13}.

\vspace{4mm}

{\bf Lemma 2.8.} \emph{Let}
$$A=\left[\begin{array}{ccccc} 0  & A_1 & & & \\ & 0  & \ddots & & \\ & & \ddots & A_{k-1} & \\ & & & 0  & B\\ & & & & C\end{array}\right] \ on \  \mathbb{C}^n=\mathbb{C}^{n_1}\oplus\cdots\oplus\mathbb{C}^{n_k}\oplus\mathbb{C}^{m},$$
\emph{where the} $A_j$'\emph{s satisfy} $A_j^*A_j=I_{n_{j+1}}$, $1\le j\le k-1$. \emph{If $W(A)$ is a circular disc centered at the origin with radius $r$ larger than} $\cos(\pi/(k+1))$, \emph{then $C$ is not invertible}.

\vspace{4mm}

{\em Proof}. Since $W(A)=\{z\in\mathbb{C}:|z|\le r\}$, $r$ is the maximum eigenvalue of $\re(e^{i\theta}A)$ and hence $\det(rI_n-\re(e^{i\theta}A))=0$ for all real $\theta$. We have
\begin{align}\label{e6}
&  0=\det\left[\begin{array}{ccccc}
rI_{n_1} & -(e^{i\theta}/2)A_1 & & & \\
-(e^{-i\theta}/2)A_1^* & rI_{n_2} & \ddots & & \\
& \ddots & \ddots & -(e^{i\theta}/2)A_{k-1} & \\
& & -(e^{-i\theta}/2)A_{k-1}^* & rI_{n_k} & -(e^{i\theta}/2)B\\
& & & -(e^{-i\theta}/2)B^* & rI_m-\re(e^{i\theta}C)\end{array}\right] \\
=& \det D_k(\theta)\cdot\det(E(\theta)-F(\theta)),\nonumber
\end{align}
where
\begin{equation}\label{e7}
D_k(\theta)=\left[\begin{array}{cccc}
rI_{n_1} & -(e^{i\theta}/2)A_1 & &  \\
-(e^{-i\theta}/2)A_1^* & rI_{n_2} & \ddots &  \\
& \ddots & \ddots & -(e^{i\theta}/2)A_{k-1}  \\
& & -(e^{-i\theta}/2)A_{k-1}^* & rI_{n_k} \end{array}\right],
  \ \  E(\theta)=rI_m-\re(e^{i\theta}C),
\end{equation}
and
\begin{equation}\label{e8}
F(\theta)=\left[0 \ \ldots \ 0 \ -(e^{-i\theta}/2)B^*\right]D_k(\theta)^{-1}\left[\begin{array}{c} 0\\ \vdots \\ 0\\ -(e^{i\theta}/2)B\end{array}\right],
\end{equation}
by using the Schur complement of $D_k(\theta)$ in the matrix in (\ref{e6}) (cf. \cite[p. 22]{8}). Note that here the invertibility of $D_k(\theta)$ follows from the facts that $D_k(\theta)$ is unitarily similar to $rI-\re J$, where $J=(\sum_{j=1}^{n_k}\oplus J_k)\oplus(\sum_{j=1}^{n_{k-1}-n_k}\oplus J_{k-1})\oplus\cdots\oplus(\sum_{j=1}^{n_1-n_2}\oplus J_1)$ (cf. the proof of Theorem 2.2 (b) $\Rightarrow$ (c)), and $r$ ($>\cos(\pi/(k+1))$) is not an eigenvalue of $\re J$. Moreover, the $(k,k)$-block of $D_k(\theta)^{-1}$ is independent of the value of $\theta$. Thus the same is true for the entries of $F(\theta)$. Under a unitary similarity, we may assume that $C=[c_{ij}]_{i,j=1}^m$ is upper triangular with $c_{ij}=0$ for all $i>j$. Let $F(\theta)=[b_{ij}]_{i,j=1}^m$ and $E(\theta)-F(\theta)=[d_{ij}(\theta)]_{i,j=1}^m$. Then
$$d_{ij}(\theta)=\left\{\begin{array}{ll}
r-\re(e^{i\theta}c_{jj})-b_{jj} \ \  \ \ & \mbox{if} \ \ i=j,\\
-(e^{i\theta}/2)c_{ij}-b_{ij}  & \mbox{if} \ \ i<j,\\
-(e^{-i\theta}/2)\overline{c}_{ji}-b_{ij} & \mbox{if} \ \ i>j.\end{array}\right.$$
Hence $p(\theta)\equiv\det(E(\theta)-F(\theta))$ is a trigonometric polynomial of degree at most $m$, say, $p(\theta)=\sum_{j=-m}^m a_je^{ij\theta}$. Since $\det(rI_m-\re(e^{i\theta}A))=0$ and $\det D_k(\theta)\neq 0$, we obtain from (\ref{e6}) that $p(\theta)=0$ for all real $\theta$. This implies that $a_j=0$ for all $j$. In particular, $a_m=(-1)^m\prod_{j=1}^m(c_{jj}/2)=0$ from the above description of the $d_{ij}(\theta)$'s. This yields that $c_{jj}=0$ for some $j$ or $C$ is not invertible. \hspace{2mm} $\blacksquare$

\vspace{4mm}

The next lemma is to be used in the proof of Lemma 2.10.

\vspace{4mm}

{\bf Lemma 2.9.} \emph{Let $A={\scriptsize\left[\begin{array}{cc} 0_p & B\\ 0 & C\end{array}\right]}$ be an $n$-by-$n$ matrix}, \emph{and let $B=[b_{ij}]_{i=1, j=1}^{p, n-p}$ and $C=[c_{ij}]_{i,j=1}^{n-p}$ with $c_{ij}=0$ for all $i>j$}. \emph{If the geometric and algebraic multiplicities of the eigenvalue $0$ of $A$ are equal to each other and $c_{11}=0$}, \emph{then $b_{i1}=0$ for all $i$}, $1\le i\le p$.

\vspace{4mm}

{\em Proof}. Let $e_j$ denote the $j$th standard unit vector $[0 \ \ldots \ 0 \ \tb{1}{j \, \mbox{th}} \ 0 \ \ldots \ 0]^T$, $1\le j\le n$. Then $e_1, \ldots, e_p$ are all in $\ker A$. Since $c_{11}=0$, we have $Ae_{p+1}=b_{11}e_1+\cdots+b_{p 1}e_p$, which is also in $\ker A$. Thus $A^2e_{p+1}=0$ or $e_{p+1}\in\ker A^2$. Our assumption on the multiplicities of 0 implies that $\ker A=\ker A^2=\cdots$. Hence we obtain $e_{p+1}\in \ker A$ or $Ae_{p+1}=0$, which yields that $b_{i 1}=0$ for all $i$, $1\le i\le p$. \hspace{2mm} $\blacksquare$

\vspace{4mm}

The following lemma is the main tool in proving, under the condition of circular $W(A)$, that $p(A)\ge a(A)-1$ yields $p(A)\ge a(A)$.

\vspace{4mm}

{\bf Lemma 2.10.} \emph{Let}
$$A=\left[\begin{array}{ccccc} 0  & A_1 & & & \\ & 0  & \ddots & & \\ & & \ddots & A_{k-2} & \\ & & & 0  & B\\ & & & & C\end{array}\right] \ \ \  \ on \  \ \mathbb{C}^n=\mathbb{C}^{n_1}\oplus\cdots\oplus\mathbb{C}^{n_{k-1}}\oplus\mathbb{C}^{m},$$
\emph{where} $k=a(A) \, (\ge 2)$, \emph{the} $A_j$'\emph{s satisfy} $A_j^*A_j=I_{n_{j+1}}$, $1\le j\le k-2$, \emph{and} $B={\scriptsize\left[\begin{array}{cc} I_p & 0\\ 0 & B_1\end{array}\right]}$ \emph{and} $C={\scriptsize\left[\begin{array}{cc} 0_p & C_1\\ 0 & C_2\end{array}\right]} \, (1\le p\le\min\{n_{k-1}, m\})$ \emph{satisfy} $B^*B+C^*C=I_m$. \emph{If $W(A)$ is a circular disc centered at the origin with radius $r$ larger than $\cos(\pi/(k+1))$}, \emph{then $A$ is unitarily similar to a matrix of the form}
$$\left[\begin{array}{ccccc} 0  & A_1' & & & \\ & 0  & \ddots & & \\ & & \ddots & A_{k-1}' & \\ & & & 0  & B'\\ & & & & C'\end{array}\right] \ \ \  \ on \  \ \mathbb{C}^n=\mathbb{C}^{n_1}\oplus\cdots\oplus\mathbb{C}^{n_{k-1}}\oplus\mathbb{C}^{q}\oplus\mathbb{C}^{m-q},$$
\emph{where} $q=\min\{n_{k-1}, m\}$, \emph{the} $A'_j$'\emph{s satisfy} ${A'_j}^*A'_j=I_{n_{j+1}}$, $1\le j\le k-2$, ${A'_{k-1}}^*A'_{k-1}=I_{q}$, \emph{and} $B'$ \emph{and} $C'$ \emph{satisfy} ${B'}^*B'+{C'}^*C'=I_{m-q}$.

\vspace{4mm}

{\em Proof}. Since $W(A)=\{z\in\mathbb{C} : |z|\le r\}$, we have $\det(rI_n-\re(e^{i\theta}A))=0$ for all real $\theta$. As in the proof of Lemma 2.8, we have the factorization $\det(rI_n-\re(e^{i\theta}A))=\det D_{k-1}(\theta)\cdot\det(E(\theta)-F(\theta))$, where $D_{k-1}(\theta)$, $E(\theta)$ and $F(\theta)$ are as in (\ref{e7}) and (\ref{e8}) with $D_{k}(\theta)^{-1}$ in the expression of $F(\theta)$ there replaced by $D_{k-1}(\theta)^{-1}$. Since $D_{k-1}(\theta)$ is unitarily similar to $rI-\re J$, where $J=(\sum_{j=1}^{n_{k-1}}\oplus J_{k-1})\oplus(\sum_{j=1}^{n_{k-2}-n_{k-1}}\oplus J_{k-2})\oplus\cdots\oplus(\sum_{j=1}^{n_1-n_{2}}\oplus J_{1})$ and the $(k-1, k-1)$-entry of $(rI_{k-1}-\re J_{k-1})^{-1}$ is $a\equiv\det(rI_{k-2}-\re J_{k-2})/\det(rI_{k-1}-\re J_{k-1})$, the $(k-1, k-1)$-block of $D_{k-1}(\theta)^{-1}$ is given by $aI_{n_{k-1}}$. Hence we have $F(\theta)=(a/4)B^*B$. As before, from $\det D_{k-1}(\theta)\neq 0$, we obtain $\det(E(\theta)-F(\theta))=0$. Thus
\begin{align}\label{e9}
& \, 0 =\det(E(\theta)-F(\theta)) \nonumber\\
= & \, \det\left(rI_m-\left[\begin{array}{cc} 0_{p} & (e^{i\theta}/2)C_1\\ (e^{-i\theta}/2)C_1^* & \re(e^{i\theta}C_2)\end{array}\right]-\frac{a}{4}\left[\begin{array}{cc} I_{p} & 0\\ 0 & B_1^*B_1\end{array}\right]\right) \nonumber\\
= & \, \det\left[\begin{array}{cc} (r-(a/4))I_p & -(e^{i\theta}/2)C_1\\ -(e^{-i\theta}/2)C_1^* & rI_{m-p}-\re(e^{i\theta}C_2)-(a/4)B_1^*B_1\end{array}\right].
\end{align}
We claim that $r\neq a/4$. Indeed, since $\det(rI_k-\re J_k)=r\det(rI_{k-1}-\re J_{k-1})-(1/4)\det(rI_{k-2}-\re J_{k-2})$, we have $\det(rI_k-\re J_k)/\det(rI_{k-1}-\re J_{k-1})=r-(a/4)$. Therefore, $r=a/4$ if and only if $\det(rI_k-\re J_k)=0$. The latter would imply $r\le\cos(\pi/(k+1))$ contradicting our assumption that $r>\cos(\pi/(k+1))$. Hence $r\neq a/4$ as asserted. Using the Schur complement, we infer from (\ref{e9}) that
$$p(\theta)\equiv\det(rI_{m-p}-\re(e^{i\theta}C_2)-\frac{a}{4}B_1^*B_1-\frac{1}{r-(a/4)}\cdot\frac{1}{4}C_1^*C_1)=0$$
for all real $\theta$. As $p(\theta)$ is a trigonometric polynomial of degree at most $m-p$, say, $p(\theta)=\sum_{j=-(m-p)}^{m-p}a_je^{ij\theta}$, this implies that $a_j=0$ for all $j$. After a unitary similarity, we may assume that $C_2=[c_{ij}]_{i,j=1}^{m-p}$ with $c_{ij}=0$ for all $i>j$. Hence $a_{m-p}=(1/2^{m-p})c_{11}\cdots c_{m-p, m-p}=0$. Thus $c_{jj}=0$ for some $j$. We may assume that $c_{11}=0$. Note that
$$A^k=\left[\begin{array}{cccc}
0 & \cdots & 0 &   (\prod_{j=1}^{k-2}A_j)BC\\
0 & \cdots & 0 &   (\prod_{j=2}^{k-2}A_j)BC^2\\
\vdots &   & \vdots   & \vdots\\
0 & \cdots & 0 &   A_{k-2}BC^{k-2}\\
0 & \cdots & 0 &   BC^{k-1}\\
0 & \cdots & 0 &   C^{k}
\end{array}\right],$$
$$BC^j=\left[\begin{array}{cc} I_p & 0\\ 0 & B_1\end{array}\right]\left[\begin{array}{cc} 0_p & C_1C_2^{j-1}\\ 0 & C_2^j\end{array}\right]=\left[\begin{array}{cc} 0_p & C_1C_2^{j-1}\\ 0 & B_1C_2^j\end{array}\right], \ \ 1\le j\le k-1,$$
and
$$C^k=\left[\begin{array}{cc} 0_p & C_1C_2^{k-1}\\ 0 & C_2^k\end{array}\right].$$
Since the first column of $C_2$ is zero, the same is true for the $(p+1)$st columns of $(\prod_{j=t}^{k-2}A_j)BC^t$ ($2\le t\le k-2$), $BC^{k-1}$ and $C^k$. As for $(\prod_{j=1}^{k-2}A_j)BC$, we need Lemma 2.9. Because $k=a(A)$, the geometric and algebraic multiplicities of the eigenvalue 0 of $A^k$ coincide. Hence we may apply Lemma 2.9 to $A^k$ to infer that the $((\sum_{j=1}^{k-1}n_j)+p+1)$st column of $A^k$ is zero. In particular, since $\ker(\prod_{j=1}^{k-2}A_j)=\{0\}$, the $(p+1)$st column of $BC={\scriptsize\left[\begin{array}{cc} 0_p & C_1\\ 0 & B_1C_2\end{array}\right]}$ is zero and thus the first column of $C_1$ is zero. Together with the zero first column of $C_2$, this yields $C={\scriptsize\left[\begin{array}{cc} 0_{p+1} & C_1^{(1)}\\ 0 & C_2^{(1)}\end{array}\right]}$. As
\begin{align*}
& I_m=B^*B+C^*C=\left[\begin{array}{cc} I_p & 0\\ 0 & B_1^*\end{array}\right]\left[\begin{array}{cc} I_p & 0\\ 0 & B_1\end{array}\right]+\left[\begin{array}{cc} 0_{p+1} & 0\\ C^{(1)*}_1 & C^{(1)*}_2\end{array}\right]\left[\begin{array}{cc} 0_{p+1} & C^{(1)}_1\\ 0 & C^{(1)}_2\end{array}\right]\\
=& \left[\begin{array}{cc} I_p & 0\\ 0 & B_1^*B_1\end{array}\right]+\left[\begin{array}{cc} 0_{p+1} & 0\\ 0 & C^{(1)*}_1C^{(1)}_1+ C^{(1)*}_2C^{(1)}_2\end{array}\right],
\end{align*}
we infer that the first column of $B_1$ is a unit vector. After another unitary similarity, we may further assume that
$$B_1=\left[\begin{array}{cc} 1 & 0\\ 0 & B_1^{(1)}\end{array}\right] \ \ \ \mbox{or} \ \ \ B=\left[\begin{array}{cc} I_{p+1} & 0\\ 0 & B^{(1)}_1\end{array}\right].$$
Applying the above arguments again, we have
$$C=\left[\begin{array}{cc} 0_{p+2} & C_1^{(2)}\\ 0 & C_2^{(2)}\end{array}\right], \ B_1^{(1)}=\left[\begin{array}{cc} 1 & 0\\ 0 & B_1^{(2)}\end{array}\right] \ \  \mbox{and} \  \ B=\left[\begin{array}{cc} I_{p+2} & 0\\ 0 & B_1^{(2)}\end{array}\right].$$
Continuing this process, we obtain
$$\mbox{(i)} \ \ C=\left[\begin{array}{cc} 0_{n_{k-1}} & C'_1\\ 0 & C'_2\end{array}\right]  \ \  \mbox{and} \  \ B=\left[ I_{n_{k-1}} \ \ 0\right] \ \ \mbox{if} \ \ n_{k-1}<m,$$
and
$$\hspace*{-25mm}\mbox{(ii)} \ \ C=0_m  \ \  \mbox{and} \  \ B=\left[\begin{array}{c} I_m \\ 0\end{array}\right] \  \ \mbox{if} \ \ n_{k-1}\ge m.$$
Finally, let $A_j'=A_j$ for $1\le j\le k-2$. In case (i), let $A_{k-1}'=I_{n_{k-1}}$, $B'=C'_1$ and $C'=C'_2$. Since
$$I_m=B^*B+C^*C=\left[\begin{array}{cc} I_{n_{k-1}} & 0\\ 0 & {C'_1}^*C'_1+{C'_2}^*C'_2\end{array}\right],$$
we have $B'^*B'+C'^*C'={C'_1}^*C'_1+{C'_2}^*C'_2=I_{m-n_{k-1}}$. On the other hand, for case (ii), let $A'_{k-1}={\scriptsize\left[\begin{array}{c} I_m \\ 0\end{array}\right]}$. In this case, $B'$ and $C'$ are absent. \hspace{2mm} $\blacksquare$

\vspace{4mm}

A consequence of the previous results is the following.

\vspace{4mm}

{\bf Proposition 2.11.} \emph{If $A$ is an $n$-by-$n$ matrix with $W(A)$ a circular disc centered at the origin and $p(A)\ge a(A)-1$}, \emph{then $p(A)=a(A)$ or $\infty$}.

\vspace{4mm}

{\em Proof}. Let $k=a(A)$. The assumption $p(A)\ge a(A)-1$ says that $A, A^2, \ldots, A^{k-1}$ are all partial isometries. In particular, we have $A^{k-1}=0$ or $\|A^{k-1}\|=1$. In the former case, $p(A)$ equals $\infty$. Hence we may assume that $\|A^{k-1}\|=1$ and thus also $\|A\|=1$. By \cite[Theorem 2.10]{1}, we have $w(A)\ge\cos(\pi/(k+1))$. Two cases are considered separately:

\vspace{4mm}

(i) $w(A)=\cos(\pi/(k+1))$. In this case, \cite[Theorem 2.10]{1} yields that $A$ is unitarily similar to a matrix of the form $J_k\oplus A_1$ with $\|A_1\|\le 1$ and $w(A_1)\le\cos(\pi/(k+1))$. Since $A_1^{k-1}$ is also a partial isometry, we may assume as before that $\|A_1^{k-1}\|=1$ and thus also $\|A_1\|=1$. Now applying \cite[Theorem 2.10]{1} again to $A_1$ yields that $w(A_1)=\cos(\pi/(k+1))$ and $A_1$ is unitarily similar to $J_k\oplus A_2$ with $\|A_2\|\le 1$ and $w(A_2)\le\cos(\pi/(k+1))$. Continuing this process, we obtain that either $p(A)=\infty$ or $A$ is unitarily similar to a direct sum of copies of $J_k$. In the latter case, we again have $p(A)=\infty$.

\vspace{4mm}

(ii) $w(A)>\cos(\pi/(k+1))$. Since $A, A^2, \ldots, A^{k-1}$ are partial isometries, Theorem 2.2 yields the unitary similarity of $A$ to a matrix of the form
$$\left[\begin{array}{ccccc} 0 & A_1 & & & \\ & 0 & \ddots & & \\ & & \ddots & A_{k-2} & \\ & & & 0 & B\\ & & & & C\end{array}\right] \ \ \  \mbox{on} \  \ \mathbb{C}^n=\mathbb{C}^{n_1}\oplus\cdots\oplus\mathbb{C}^{n_{k-1}}\oplus\mathbb{C}^{m}
$$
with $A_j^*A_j=I_{n_{j+1}}$, $1\le j\le k-2$, and $B^*B+C^*C=I_m$. By Lemma 2.8, $C$ is not invertible. We may assume, after a unitary similarity, that $B$ and $C$ are of the forms ${\scriptsize\left[\begin{array}{cc} 1 & 0\\ 0 & B_1 \end{array}\right]}$ and ${\scriptsize\left[\begin{array}{cc} 0 & C_1\\ 0 & C_2\end{array}\right]}$, where $B_1$, $C_1$ and $C_2$ are $(n_{k-1}-1)$-by-$(m-1)$, 1-by-$(m-1)$ and $(m-1)$-by-$(m-1)$ matrices, respectively. Using Lemma 2.10, we obtain the unitary similarity of $A$ to a matrix of the form in Theorem 2.2 (b). Thus, by Theorem 2.2 again, $A, A^2, \ldots, A^k$ are partial isometries. Hence $p(A)\ge k=a(A)$. Our assertion then follows from Corollary 2.5. \hspace{2mm} $\blacksquare$

\vspace{4mm}

Note that, in the preceding proposition, the number ``$a(A)-1$'' is sharp as was seen from Example 2.7.

\vspace{4mm}

We are now ready to prove Theorem 2.6.

\vspace{4mm}

{\em Proof of Theorem $2.6$}. The implications (b) $\Rightarrow$ (c) and (c) $\Rightarrow$ (d) are trivial. On the other hand, (d) $\Rightarrow$ (a) follows from Theorem 2.4. Hence we need only prove (a) $\Rightarrow$ (b). Let $k=a(A)$. By Proposition 2.11, $A, A^2, \ldots, A^k$ are partial isometries. Thus $A$ is unitarily similar to the matrix $A'$ in Theorem 2.2 (b). Since $k$ is the ascent of $A$, the geometric multiplicity of $A^k$, that is, $\nul A^k$ is equal to the algebraic multiplicity of eigenvalue 0 of $A$. As proven in (a) $\Rightarrow$ (b) of Theorem 2.2, $\nul A^k=\sum_{j=1}^kn_j$. We infer from the structure of $A'$ that 0 cannot be an eigenvalue of $C$. On the other hand, applying Lemma 2.8 to $A'$ yields the noninvertibility of $C$. This leads to a contradiction. Thus $B$ and $C$ won't appear in $A'$ and, therefore, $A'$, together with $A$, is unitarily similar to a direct sum of Jordan blocks by Theorem 2.2 (c). This proves (b). \hspace{2mm} $\blacksquare$

\vspace{8mm}

\noindent
{\bf\large 3. $S_n$-matrices}

\vspace{4mm}

In this section, we apply the results in Section 2 to the class of $S_n$-matrices. This we start with the following.

\vspace{4mm}

{\bf Proposition 3.1.} \emph{Let $A$ be a noninvertible $S_n$-matrix}. \emph{Then}

(a) \emph{$a(A)$ equals the algebraic multiplicity of the eigenvalue $0$ of $A$},

(b) $p(A)=a(A)$ \emph{or} $\infty$,

(c) \emph{$p(A)=\infty$ if and only if $A$ is unitarily similar to $J_n$}, \emph{and}

(d) $\rank A^j=n-j$ \emph{for} $1\le j\le a(A)$.

\vspace{4mm}

{\em Proof}. Let $k= a(A)$.

\vspace{4mm}

(a) It is known that, for any eigenvalue $\lambda$ of $A$, there is exactly one associated block, say, $\lambda I_{\ell}+J_{\ell}$ in the Jordan form of $A$. In particular, for $\lambda=0$, both $a(A)$ and the algebraic multiplicity of 0 are equal to the size $\ell$ of its associated Jordan block $J_{\ell}$.

\vspace{4mm}

(b) By \cite[Corollary 1.3]{2}, $A$ is unitarily similar to a matrix of the form $A'\equiv{\scriptsize\left[\begin{array}{cc} J_k & B\\ 0 & C\end{array}\right]}$, where $B={\scriptsize\left[\begin{array}{c} 0\\ b\end{array}\right]}$ is a $k$-by-$(n-k)$ matrix with $b$ a row vector of $n-k$ components, and $C$ is an invertible $(n-k)$-by-$(n-k)$ upper-triangular matrix. Since $\rank(I_n-A^*A)=1$, we infer from
$$I_n-A'^*A'=\left[\begin{array}{cc} I_k & 0\\ 0 & I_{n-k}\end{array}\right]-\left[\begin{array}{cc} J_k^* & 0\\ B^* & C^*\end{array}\right]\left[\begin{array}{cc} J_k & B\\ 0 & C\end{array}\right]=\left[\begin{array}{cc} {\scriptsize\left[\begin{array}{cccc} 1 & & &\\ & 0 & & \\ & & \ddots & \\ & & & 0\end{array}\right]} & 0\\ 0 & I_{n-k}-(B^*B+C^*C)\end{array}\right]$$
that $B^*B+C^*C=I_{n-k}$. As $A'$ can also be expressed as
$$\left[\begin{array}{ccccc}
0 & 1 &   &  & 0\\
  & 0 & \ddots &   & \vdots\\
  &   & \ddots & 1 & 0\\
  &   &   &  0 & b  \\
  &   & & &  C
\end{array}\right] \ \ \ \mbox{on} \ \ \mathbb{C}^n=\underbrace{\mathbb{C}\oplus\cdots\oplus \mathbb{C}}_k\oplus \mathbb{C}^{n-k}$$
with $b^*b+C^*C=I_{n-k}$, Theorem 2.2 can be invoked to conclude that $A, A^2, \ldots, A^k$ are partial isometries. Thus $p(A)\ge k$. It follows from Corollary 2.5 that $p(A)=k$ or $\infty$.

\vspace{4mm}

(c) If $p(A)=\infty$, then the unitary similarity of $A$ and $J_n$ is an easy consequence of Theorem 2.4 and the fact that $A$ is irreducible (in the sense that it is not unitarily similar to the direct sum of two other matrices). The converse is trivial.

\vspace{4mm}

(d) As in the proof of (b), $A$ is unitarily similar to $A'=\left[\begin{array}{cc} J_k & B\\ 0 & C\end{array}\right]$, where $B={\scriptsize\left[\begin{array}{c} 0\\ b\end{array}\right]}$ and $C$ is invertible. Then $A^j$ is unitarily similar to
$$A'^j=\begin{array}{ll} \ \ \ \overbrace{\ \hspace{15mm} \ }^{\displaystyle j} \  \ \overbrace{\ \hspace{23mm} \ }^{\displaystyle k-j} & \\ \left[\begin{array}{c|c}
       \begin{array}{ccccccc}
         0 & \cdots & 0 & 1 & 0 & \cdots & 0   \\
           & \cdot  &   & 0 & \ddots & \ddots & \vdots   \\
           &   & \cdot  &   & \ddots  & \ddots & 0   \\
           &   &   & \cdot  &   & \ddots  & 1  \\
           &   &   &   & \cdot  &   & 0  \\
           &   &   &   &   & \cdot  & \vdots  \\
           &   &   &   &   &   & 0
       \end{array} & \begin{array}{c} \\ 0  \vspace{3mm} \\ \\ \\ B_j\end{array}\\ \hline 0 & C^j\end{array}
     \right] & \hspace{-11mm}\begin{array}{l}\vspace{2mm}\left.\begin{array}{l} {\ } \\ {\ } \\ {\ } \\ {\ }\end{array}\right\}k-j\\ \vspace{5mm}\left.\begin{array}{l}{\ } \\ {\ } \\ \vspace*{-2mm}{\ }\end{array}\right\}j\end{array}\end{array}.$$
for some $j$-by-$(n-k)$ matrix $B_j$. Since the first $k-j$ rows and the last $n-k$ rows of $A'^j$ are linearly independent, we infer that $\rank A^j=\rank A'^j=(k-j)+(n-k)=n-j$ for $1\le j\le k$. \hspace{2mm} $\blacksquare$

\vspace{4mm}

The next corollary complements Corollary 2.5: it shows that any allowable value for $p(A)$ can actually be attained by some matrix $A$.

\vspace{4mm}

{\bf Corollary 3.2.} \emph{For any integers $n$ and $j$ satisfying $1\le j\le n-1$}, \emph{there is an $n$-by-$n$ matrix $A$ with $p(A)=j$}.

\vspace{4mm}

{\em Proof}. Let $A$ be a noninvertible $S_n$-matrix with the algebraic multiplicity of its eigenvalue 0 equal to $j$ (cf. \cite[Corollary 1.3]{2}). Then $p(A)=a(A)=j$ by Proposition 3.1. \hspace{2mm} $\blacksquare$

\vspace{4mm}

For an $n$-by-$n$ matrix $A=[a_{ij}]_{i,j=1}^n$ and an $m$-by-$m$ matrix $B$, their {\em tensor product} (or {\em Kronecker product}) $A\otimes B$ is the $(nm)$-by-$(nm)$ matrix
$$\left[
    \begin{array}{ccc}
      a_{11}B & \cdots & a_{1n}B \\
      \vdots &   & \vdots \\
      a_{n1}B & \cdots & a_{nn}B
    \end{array}
  \right].$$
Basic properties of tensor products can be found in \cite[Chapter 4]{9}. Our main concern here is when $W(A)$ and $W(A\otimes A)$ are circular discs (centered at the origin). Problems of this nature have also been considered in \cite{1}. The main result of this section is the following theorem.

\vspace{4mm}

{\bf Theorem 3.3.} \emph{Let $A$ be an $S_n$-matrix}. \emph{Then the following conditions are equivalent}:

(a) \emph{$W(A)$ is a circular disc centered at the origin},

(b) \emph{$W(A\otimes A)$ is a circular disc centered at the origin}, \emph{and}

(c) \emph{$A$ is unitarily similar to $J_n$}.

\vspace{4mm}

In preparation for its proof, we need the next lemma.

\vspace{4mm}

{\bf Lemma 3.4.} \emph{Let $A$ and $B$ be an $n$-by-$n$ and $m$-by-$m$ nonzero matrices}, \emph{respectively}.

(a) $$a(A\otimes B)=\left\{\begin{array}{ll}
\min\{a(A), a(B)\} \ \ \ & \mbox{\em if} \ \ a(A), a(B)\ge 1,\\
a(A) & \mbox{\em if} \ \ a(B)=0,\\
a(B) & \mbox{\em if} \ \ a(A)=0.\end{array}\right.$$

(b) \emph{If $A$ and $B$ are partial isometries}, \emph{then so is $A\otimes B$}. \emph{The converse is false}.

(c) \emph{Assume that $A$ and $B$ are} (\emph{nonzero}) \emph{contractions}. \emph{Then $A$ and $B$ are partial isometries if and only if $A\otimes B$ is a partial isometry}.

(d) \emph{If $A$ and $B$ are} (\emph{nonzero}) \emph{contractions}, \emph{then} $p(A\otimes B)=\min\{p(A), p(B)\}$.

(e) \emph{$A$ is a partial isometry if and only if $A\otimes A$ is}. \emph{Thus}, \emph{in particular}, $p(A\otimes A)=p(A)$.

\vspace{4mm}

The proof makes use of the facts that (i) if $A$ (resp., $B$) is similar to $A'$ (resp., $B'$), then $A\otimes B$ is similar to $A'\otimes B'$, and (ii) if the eigenvalues of $A$ (resp., $B$) are $a_i$, $1\le i\le n$ (resp., $b_j$, $1\le j\le m$), then the eigenvalues of $A\otimes B$ are $a_ib_j$, $1\le i\le n$, $1\le j\le m$, counting algebraic multiplicities (cf. \cite[Theorem 4.2.12]{9}).

\vspace{4mm}

{\em Proof of Lemma $3.4$}. (a) Let $k_1=a(A)$ and $k_2=a(B)$, and assume that $2\le k_1\le k_2$. Let $J_{k_1}$ (resp., $J_{k_2}$) be a Jordan block in the Jordan form of $A$ (resp., $B$). Since
$$(J_{k_1}\otimes J_{k_2})^{k_1}=J_{k_1}^{k_1}\otimes J_{k_2}^{k_1}=0_{k_1}\otimes J_{k_2}^{k_1}=0_{k_1k_2}$$
and
$$(J_{k_1}\otimes J_{k_2})^{k_1-1}=J_{k_1}^{k_1-1}\otimes J_{k_2}^{k_1-1}\neq 0_{k_1k_2},$$
the size of the largest Jordan block in the Jordan form of $A\otimes B$ is $k_1$. This shows that $a(A\otimes B)=k_1=\min\{a(A), a(B)\}$. The other cases can be proven even easier.

\vspace{4mm}

(b) This is a consequence of the equivalence of (a) and (b) in Lemma 2.1 as $A^*A$ and $B^*B$ are projections, which implies the same for $(A\otimes B)^*(A\otimes B)$. The converse is false as seen by the example of $A=[2]$ and $B=[1/2]$.

\vspace{4mm}

(c) If $A\otimes B$ is a partial isometry, then $(A\otimes B)^*(A\otimes B)=(A^*A)\otimes(B^*B)$ is a projection by Lemma 2.1. Since the positive semidefinite $A^*A$ and $B^*B$ are both contractions, their eigenvalues $a_i$, $1\le i\le n$, and $b_j$, $1\le j\le m$, are such that $0\le a_i, b_j\le 1$ for all $i$ and $j$. As the eigenvalues of $(A^*A)\otimes(B^*B)$, the products $a_ib_j$, $1\le i\le n$, $1\le j\le m$, can only be $0$ and $1$. Thus the same is true for the $a_i$'s and $b_j$'s. It follows that $A^*A$ and $B^*B$ are projections. Therefore, $A$ and $B$ are partial isometries.

\vspace{4mm}

(d) This follows from (c) immediately.

\vspace{4mm}

(e) If $A\otimes A$ is a partial isometry, then $(A\otimes A)^*(A\otimes A) = (A^*A)\otimes (A^*A)$ is a projection with eigenvalues 0 and 1.  But its eigenvalues are also given by $a_i a_j$, $1\le i, j \le n$, where the $a_i$'s are eigenvalues of $A^*A$.  If any $a_i$ is nonzero and not equal to 1, then the same is true for $a_i^2$, which is a contradiction.  Hence all the $a_i$'s are either 0 or 1.  It follows that $A^*A$ is a projection and $A$ is a partial isometry.  The converse was proven in (c). \hspace{2mm} $\blacksquare$

\vspace{4mm}

Finally, we are ready to prove Theorem 3.3.

\vspace{4mm}

{\em Proof of Theorem $3.3$}. To prove (a) $\Rightarrow$ (c) (resp., (b) $\Rightarrow$ (c)), note that the center of the circular $W(A)$ (resp., $W(A\otimes A)$) must be an eigenvalue of $A$ (resp., $A\otimes A$) (cf. \cite[Theorem]{3}). In particular, this says that $A$ (resp., $A\otimes A$) is noninvertible. Since the eigenvalues of $A\otimes A$ are $a_ia_j$, $1\le i, j\le n$, where the $a_i$'s are the eigenvalues of $A$ (cf. \cite[Theorem 4.2.12]{9}), the noninvertibility of $A\otimes A$ also implies that of $A$. Hence $p(A)=a(A)$ or $\infty$ by Proposition 3.1 (b). If $p(A)=\infty$, then we have already had (c) by Proposition 3.1 (c). Thus we may assume that $p(A)=a(A)$. In this case, we also have
$$p(A\otimes A)=p(A)=a(A)=a(A\otimes A)$$
by Lemma 3.4 (d) (or (e)) and (a). Applying Theorem 2.6, we obtain the unitary similarity of $A$ (resp., $A\otimes A$) to a direct sum of Jordan blocks. It follows that the only eigenvalue of $A$ (resp., $A\otimes A$ and hence of $A$) is 0. Hence $A$ is unitarily similar to $J_n$, that is, (c) holds.

\vspace{4mm}

The implication (c) $\Rightarrow$ (a) is trivial since, under (c), we have $W(A)=\{z\in\mathbb{C} : |z|\le\cos(\pi/(n+1))\}$. For (c) $\Rightarrow$ (b), note that (c) implies that $A$ is unitarily similar to $e^{i\theta}A$ for all real $\theta$. Hence $A\otimes A$ is unitarily similar to $e^{i\theta}(A\otimes A)$ for real $\theta$. Thus $W(A\otimes A)$ is a circular disc centered at the origin. This also follows from \cite[Proposition 2.8]{1}. \hspace{2mm} $\blacksquare$

\vspace{4mm}

We remark that the equivalence of (a) and (c) in Theorem 3.3 was shown before in \cite[Lemma 5]{12} by a completely different proof.

\vspace{4mm}

We end this section with two examples and one open question. The examples show that, in contrast to the case of $S_n$-matrices, the conditions of $W(A)$ and $W(A\otimes A)$ being circular discs centered at the origin are independent of each other for a general matrix $A$.

\vspace{4mm}

{\bf Example 3.5.} Let $A=[\lambda]\oplus J_2$, where $1/2<|\lambda|\le 1/\sqrt{2}$. Then
$$W(A\otimes A)=W([\lambda^2]\oplus\lambda J_2\oplus\lambda J_2\oplus\left[\begin{array}{cc} 0_2 & J_2\\ 0 & 0_2\end{array}\right])=\{z\in\mathbb{C} : |z|\le\frac{1}{2}\},$$
but $W(A)$, being the convex hull of $\{\lambda\}\cup\{z\in\mathbb{C} : |z|\le 1/2\}$, is obviously not a circular disc.

\vspace{4mm}

{\bf Example 3.6.} Let $$A=\left[\begin{array}{ccc} 0 & -\sqrt{2} & 1\\ 0 & 0 & 1\\ 0 & 0 & \sqrt{2}/2\end{array}\right].$$
Then, for any real $\theta$,
$$\re(e^{i\theta}A)=\frac{1}{2}\left[\begin{array}{ccc} 0 & -\sqrt{2}e^{i\theta} & e^{i\theta}\\ -\sqrt{2}e^{-i\theta} & 0 & e^{i\theta}\\ e^{-i\theta} & e^{-i\theta} & \sqrt{2}\cos\theta\end{array}\right],$$
whose maximum eigenvalue can be computed to be always equal to 1. Hence $W(A)=\overline{\mathbb{D}}$. On the other hand, a long and tedious computation shows that the characteristic polynomial $p(z)\equiv\det(zI_9-2\re(A\otimes A))$ of $2\re(A\otimes A)$ can be factored as
\begin{equation}\label{e10}
z^2(z^2-3)(z^5-z^4-17z^3+17z^2+46z-48).
\end{equation}
Assume that $W(A\otimes A)=\{z\in\mathbb{C} : |z|\le \sqrt{r}/2\}$ for some $r>0$. Then the maximum and minimum eigenvalues of $2\re(A\otimes A)$ are $\sqrt{r}$ and $-\sqrt{r}$, respectively. Note that $p(2)=-8<0$ and $p(\infty)=\infty$ imply that $p$ has a zero larger than 2. Hence $r\neq 3$. Similarly, we have $r\neq -3$. Since both $\sqrt{r}$ and $-\sqrt{r}$ are zeros of $p$, we also have
\begin{align}\label{e11}
& \ p(z)= z^2(z^2-3)(z^2-r)(z^3+az^2+bz+c) \nonumber\\
=& \ z^2(z^2-3)(z^5+az^4+(b-r)z^3+(c-ar)z^2-brz-cr)
\end{align}
for some real $a$, $b$ and $c$. Comparing the coefficients of the last factors in (\ref{e10}) and (\ref{e11}) yields that $a=-1$, $b-r=-17$, $c-ar=17$, $br=-46$ and $cr=48$. From these, we deduce that $c+r=17$ and hence $b=-c$. This leads to $-46=br=-cr$, which contradicts $cr=48$. Thus $W(A\otimes A)$ cannot be a circular disc at 0.

\vspace{4mm}

The matrix $A$ in the preceding example was also considered in \cite[Example 3.4]{1} for another purpose.

\vspace{4mm}

{\bf Question 3.7.} Is it true that, for any integers $n$, $j$ and $k$ satisfying $1\le j\le k\le n-1$, there is an $n$-by-$n$ matrix $A$ with $p(A)=j$ and $a(A)=k$? This is a refinement of Corollary 3.2. It is true if $k<n/2$. Indeed, in this case, we have $j\le k\le n-k-1$. Let $A=J_k\oplus B$, where $B$ is a noninvertible $S_{n-k}$-matrix whose eigenvalue 0 has algebraic multiplicity $j$. Then $p(A)=p(B)=a(B)=j$ by Proposition 3.1. On the other hand, we obviously have $a(A)=k$.

\end{document}